\documentclass{amsart}
\usepackage{amsmath,amsthm,amssymb,amsfonts}

\makeatletter
\def\proof{\@ifstar{P\,r\,o\,o\,f}{P\,r\,o\,o\,f.\ }}

\renewcommand\th@remark{%
  \thm@headfont{\bfseries}%
  \normalfont % body font
  \thm@preskip\topsep \divide\thm@preskip\tw@
  \thm@postskip\thm@preskip
}

\renewenvironment{equation}{\refstepcounter{equation}$$}{\eqno{(\thesection.\theequation)}$$}

\makeatother

\newcounter{Example}[section]
\newcounter{Th}[section]
\newcounter{Pr}[section]
\newcounter{Lm}[section]
\newcounter{Remark}[section]
\newcounter{Def}[section]
\newcounter{lcounter}[section]
\newcounter{Corol}[section]

\newenvironment{Th}[1][\relax]
    {\medspace\refstepcounter{Th}T\,h\,e\,o\,r\,e\,m \arabic{section}.\theTh.\ \it}
    {\rm\medspace}

\newenvironment{Th.}[1][\relax]
    {\medspace\refstepcounter{Th}T\,h\,e\,o\,r\,e\,m \arabic{section}.\theTh.\ \it}
    {\rm\medspace}

\newenvironment{Pr.}[1][\relax]
    {\medspace\refstepcounter{Pr}P\,r\,o\,p\,o\,s\,i\,t\,i\,o\,n \arabic{section}.\thePr.\ \it}
    {\rm\medspace}

\newenvironment{Lm}[1][\relax]
    {\medspace\refstepcounter{Lm}L\,e\,m\,m\,a \arabic{section}.\theLm.\ \it}
    {\rm\medspace}

\newenvironment{Corol}[1][\relax]
    {\medspace\refstepcounter{Corol} C\,o\,r\,o\,l\,l\,a\,r\,y \arabic{section}.\theCorol.\ \it}
    {\rm\medspace}

\newenvironment{Remark}[1][\relax]
    {\medspace\refstepcounter{Remark}R\,e\,m\,a\,r\,k \arabic{section}.\theRemark.\rm\ }
    {\medspace}
\newenvironment{Example}[1][\relax]
    {\medspace\refstepcounter{Example}E\,x\,a\,m\,p\,l\,e \arabic{section}.\theExample.\rm\ }
    {\medspace}

\newenvironment{Def}[1][\relax]
    {\medspace\refstepcounter{Def}D\,e\,f\,i\,n\,i\,t\,i\,o\,n \arabic{section}.\theDef.\rm\ }
    {\medspace}

\newenvironment{Def.}[1][\relax]
    {\medspace\refstepcounter{Def}D\,e\,f\,i\,n\,i\,t\,i\,o\,n \arabic{section}.\theDef.\rm\ }
    {\medspace}

\def\R{{\mathbb R}}  % множество вещественных чисел
\def\H{\mathcal{H}} % гильбертово пространство
\def\cl{\mathop{\rm cl\,}}

\def\conv{\mathop{\rm conv\,}}

\def\l{\langle}
\def\r{\rangle}
\def\ll{\lambda}
\def\d {\partial\,}
\def\ep{\varepsilon} % эпсилон
\def\Int {\mbox{\rm int\,}} % внутренность

\def\L {{\mathcal L}}  % каллиграфическое L
\def\diam {\mbox{\rm diam}\,}

\begin{document}

\title[Uniform convexity and the splitting problem for
selections]{Uniform convexity and the splitting problem for
 selections}

\author[M. V. Balashov and D. Repov\v{s}]{Maxim V. Balashov and Du\v{s}an Repov\v{s}}

%\thanks{Supported by }

\address{Department of Higher Mathematics, Moscow Institute of Physics and Technology, Institutski str. 9,
Dolgoprudny, Moscow region, Russia 141700. balashov@mail.mipt.ru}
\address{Faculty of Mathematics and Physics, and Faculty of Education, University of Ljubljana, Jadranska 19, Ljubljana, Slovenia 1000.
dusan.repovs@guest.arnes.si}



\begin{abstract}
We continue to investigate cases when the Repov\v{s}-Semenov
splitting problem for selections has an
affirmative solution for
continuous set-valued mappings. We consider the situation in
infinite-dimensional uniformly convex Banach spaces. We use the
notion of Polyak of uniform convexity and modulus of
uniform convexity for  arbitrary convex sets (not necessary
balls).  We study general geometric properties of uniformly
convex sets. We also obtain an affirmative
solution of the splitting
problem for selections of certain
set-valued mappings with uniformly
convex images.

\end{abstract}
\date{\today}
\maketitle

\section{Introduction}

\def\i {\mbox{\rm int}\,}

The questions concerning continuity of set-valued mappings and
existence of continuous, uniformly continuous
and Lipschitz continuous selections of set-valued mappings
have for a long time been
the central questions of nonsmooth analysis
\cite{Aubin}, \cite{Aubin+Frankowska}.
The classical Michael theorem \cite{MichaelII} guarantees
the existence
of continuous selections for lower semicontinuous set-valued
mappings with convex closed images.
However, the condition of lower
semicontinuity for
a set-valued mapping is not typical for (many)
problems in which the set-valued mappings are
represented
as the intersection of two set-valued
mappings. This occurs e.g.
in approximation theory
\cite{Berd}, \cite{Mar}.

It is well-known (\cite{Aubin}, \cite{Aubin+Frankowska}) that
even the intersection of Lipschitz continuous set-valued mappings
with convex compact images, defined on  $\R^{n}$,
is only upper
semicontinuous.
In certain minimization problems \cite{Polyak} and
problems of stability of functionals \cite{Berd} it is necessary
to obtain uniformly continuous selections and explicit
estimates for their moduli of continuity.
This explains the necessity for additional constraints on
the type of convexity of the set-valued mappings
under consideration.

Let $E$ be a Banach space. The {\it diameter}
of the subset $A\subset E$
is defined as $\diam A=\sup\limits_{x_{1},x_{2}\in A}\|
x_{1}-x_{2}\|$. Let $\d A$ be the boundary of the set $A$, $\i A$
the interior of $A$, and $\cl A$ the closure of $A$. Let $\l p,x
\r$ be the value of the functional $p\in E^{*}$ at the point $x\in
E$. We define the {\it closed ball} with center $a\in E$ and radius $r$
as follows: $B_{r}(a)=\{ x\in E\ |\ \| x-a\|\le r\}$.
Following \cite{Polyak}, we define uniformly
convex set as follows:

\begin{Def}\label{modulus} (\cite{Polyak})
Let $E$ be a Banach space and  $A\subset E$ a closed convex set.
\it The modulus of convexity  \rm $\delta_{A}:\ (0,\diam A)\to
[0,+\infty) $ is the function defined by
$$
\delta_{A}(\ep) = \sup\left\{ \delta\ge 0\ \left|\
B_{\delta}\left( \frac{x_{1}+x_{2}}{2}\right)\right.\subset A,\
\forall x_{1},x_{2}\in A:\ \| x_{1}-x_{2}\|=\ep \right\}.
$$\rm
\end{Def}
\begin{Def}\label{RM} (\cite{Polyak})
Let $E$ be a Banach space and  $A\subset E$ a closed convex set.
If the modulus of convexity  $\delta_{A}(\ep )$ is strictly
positive for all $\ep\in (0,\diam A)$, then we call the set $A$
\it uniformly convex (with modulus $\delta_{A}(\cdot)$).\rm
\end{Def}

Definition 1.\ref{modulus} is very similar to the well-known
definition of the modulus of convexity for uniformly convex
function \cite[Chapter 4 \S 7]{Vas}. If the set $A$ is bounded
and has the center of symmetry then $\delta_{A}(\cdot)$ is the
modulus of convexity for space $E$ with the ball $A$
(\cite{Diestel}, \cite{Lindestrauss+tzafriri}). Note that, as in
the case of the bodies with center of symmetry (under assumption
$A\ne E$, see \cite[Part e]{Lindestrauss+tzafriri}), it
suffices to choose points $x_{1},x_{2}\in \d A$, i.e.
$$
\delta_{A}(\ep) = \sup\left\{ \delta\ge 0\ \left|\
B_{\delta}\left( \frac{x_{1}+x_{2}}{2}\right)\right.\subset A,\
\forall x_{1},x_{2}\in \d A:\ \| x_{1}-x_{2}\|=\ep \right\}.
$$

The properties of uniformly convex sets were used in \cite{PL, Polyak}
for the proof of convergence of minimizing sequences in certain
extremal problems. Similar constructions appeared in
approximation theory (see for example \cite{Berd}, \cite[P.
12]{Berd-d}).
We plan to consider the entire class of uniformly convex sets and
apply their properties for the solution of the splitting problem
for selections.

The splitting problem for selections was formulated in
\cite{Repovs+Semenov}. Let $F_{i}:X\to 2^{Y_{i}}$, $i=1,2$, be
any (lower semi)continuous mappings with closed convex images and
let $L:Y_{1}\oplus Y_{2}\to Y$ be any linear surjection. The
splitting problem is the problem of representing any continuous
selection $f\in L(F_{1},F_{2})$ in the form $f= L(f_{1},f_{2})$,
where $f_{i}\in F_{i}$ are some continuous selections, $i=1,2$.
Some special cases of this problem in finite-dimensional spaces
were considered in \cite{MM}, \cite{RS}.

In \cite{Bal+Rep} we obtained new results for finite-dimensional
spaces and proved that there exist approximate solutions of the
splitting problem for Lipschitz selections in the Hilbert space.
We also wish to mention
\cite{LM-LS} and \cite{M-LS}, where related
questions were considered.

\section{Uniformly convex sets and their properties}

Note that if a set is uniformly convex then it is also strictly
convex, i.e. its boundary
contains
no nondegenerate
segments.

\begin{Lm}\label{1}
Let $A\subset E$ be a closed and uniformly convex set with modulus
$\delta_{A} (\cdot)$ and suppose that $A\ne E$. Then for any
$\lambda\in (0,1)$, $\ep\in (0,\diam A)$ the following inequality
holds
$$
\delta_{A}(\lambda \ep) \le \lambda\delta_{A} (\ep).
$$
\end{Lm}

Note that for any uniformly convex unit ball $A$, the inequality
$\delta_{A}(\lambda \ep) \le \lambda\delta_{A} (\ep)$, for all
$\lambda, \ep\in (0,1)$, follows from \cite[Lemma
1.e.8]{Lindestrauss+tzafriri}.

\proof Let's fix $\ep\in (0,\diam A)$, $\alpha>0$ and $\lambda\in
(0,1)$. Choose points $x_{1},x_{2}\in \d A$, such that $\|
x_{1}-x_{2}\|=\ep$ and $\delta_{A} (\ep)+\alpha>\delta$, where
$\delta = \sup\{ {r\ge 0}  \ |\ B_{r}(z)\subset A\}$ and
$z=\frac12 (x_{1}+x_{2})$.

For any $k$ we define a point $a_{k}\in \d A$ with $\|
a_{k}-z\|\le \delta+\frac1k$. Let $y^{k}_{i}$ be the homothetic
image of the point $x_{i}$ under the homothety with center $a_{k}$
and coefficient $\lambda$, $i=1,2$; let $z_{k}$ be the homothetic
image of the point $z$ under the homothety with center $a_{k}$ and
coefficient $\lambda$.

We have $\| y^{k}_{1}-y^{k}_{2}\|=\lambda \ep$ and $\|
z_{k}-a_{k}\|\le \lambda \delta+\lambda\frac1k$.
It follows
from the inclusions $y_{i}^{k}\in A$, $i=1,2$,  that
$$
\delta_{A}(\ll\ep)\le \| z_{k}-a_{k}\|\le
\lambda\delta+\lambda\frac{1}{k}\le\lambda\delta_{A}(\ep)+\lambda\alpha+\lambda\frac{1}{k}.
$$
By taking limits $\alpha\to+0$, $k\to\infty$ we get the following
inequality:
$$
\delta_{A}(\ll\ep)\le \lambda\delta_{A}(\ep ).
$$
\qed

The following corollary follows from Lemma 2.\ref{1}.

\begin{Corol}\label{strogo}
The modulus of convexity is a strictly monotone function and
moreover, the function $\ep \to \frac{\delta_{A}(\ep)}{\ep}$ is
also monotone.
\end{Corol}

\begin{Lm}\label{2}
Let $A\subset E$ be a closed and uniformly convex set with modulus
$\delta_{A} (\cdot)$. Let $\ep\in (0,\diam A)$, $p_{1},p_{2}\in \d
B^{*}_{1}(0)$, $x_{i}=\arg\max\limits_{x\in A}\l p_{i},x\r$,
$i=1,2$. If $\| p_{1}-p_{2}\|<\frac{4\delta_{A}(\ep )}{\ep}$ then
$\| x_{1}-x_{2}\|< \ep$.
\end{Lm}

\proof Suppose that $\| x_{1}-x_{2}\|\ge\ep$. Define $\delta =
\delta_{A} (\| x_{1}-x_{2}\|)$. We have $B_{\delta}\left(
\frac{x_{1}+x_{2}}{2}\right)\subset A$. By hypotheses of the
lemma,
$$ \l p_{1},x_{1}\r = \max\limits_{x\in A}\l p_{1},x\r\ge
\max\limits_{x\in B_{\delta}(\frac{x_{1}+x_{2}}{2})}\l p_{1},x\r\
=  \frac12 \l p_{1},x_{1}+x_{2}\r+\delta$$
and in the same way $\l
p_{2},x_{2}\r\ge \frac12 \l p_{2},x_{1}+x_{2}\r+\delta$. Hence
$$
\l p_{1},x_{1}\r-\l p_{1},x_{2}\r\ge 2\delta, \quad \l
p_{2},x_{2}\r-\l p_{2},x_{1}\r\ge 2\delta.$$
Adding the
last two inequalities
$$
\l p_{1}-p_{2},x_{1}-x_{2}\r\ge 4\delta.
$$
we obtain $\| p_{1}-p_{2}\|\cdot\| x_{1}-x_{2}\|\ge 4\delta$ and
$$
\| p_{1}-p_{2}\|\ge \frac{4\delta_{A}(\| x_{1}-x_{2}\|)}{\|
x_{1}-x_{2}\|}\ge\frac{4\delta_{A}(\ep )}{\ep},
$$
where the last
inequality follows by Corollary 2.\ref{strogo}. \qed

Let us denote $\varphi (\ep)=\frac{4\delta_{A}(\ep)}{\ep}$. We
obtain the following corollary:

\begin{Corol}\label{Holder}
Let $A\subset E$ be a closed and uniformly convex set with modulus
$\delta_{A} (\cdot)$. Let $p_{1},p_{2}\in \d B^{*}_{1}(0)$,
$x_{i}=\arg\max\limits_{x\in A}\l p_{i},x\r$, $i=1,2$. Then
$$
\varphi (\| x_{1}-x_{2}\|)\le \| p_{1}-p_{2}\|.
$$
\end{Corol}
\proof Let $\| x_{1}-x_{2}\|=\ep$. By Lemma 2.\ref{2} we then
obtain that $\varphi (\ep)=\frac{4\delta_{A}(\ep)}{\ep}\le \|
p_{1}-p_{2}\|$.\qed

\begin{Remark}\label{Midresult}
Suppose that the convex closed bounded subset $A$ of a Banach
space $E$ has uniformly continuous supporting elements,
i.e. that
there exists a continuous function $\varphi: [0,\diam A)\to
[0,+\infty)$, $\varphi (0)=0$, such that for any unit vectors
$p_{1},p_{2}\in E^{*}$ and $x_{i}=\arg\max\limits_{x\in A}\l
p_{i},x\r$, $i=1,2$:
$$
\varphi (\| x_{1}-x_{2}\|)\le \| p_{1}-p_{2}\|.
$$
Then there exists $C>0$ such that
$$
\delta_{A}(\ep)\ge C\cdot\int\limits_{0}^{\frac{\ep}{2}}\varphi
(t)\, dt,\qquad \forall\ep\in (0,\diam A).
$$
The proof of this fact has not been
published yet, however, it is too long to
be
included in this paper.
\end{Remark}

The supporting function of the set $A\subset E$ is defined by
$s(p,A)=\sup\limits_{x\in A}\l p,x\r$, $p\in E^{*}$. This is a
positively uniform convex closed function (see
\cite{Aubin,Polovinkin&Balashov}). For the set $A$ we define the
barrier cone by $b(A)=\{ p\in E^{*}\ |\ s(p,A)<+\infty\}$, i.e. $b(A)$
is the domain of the supporting function.

The fact that every
uniformly convex set which does not coincide
with the entire space is bounded was stated
in \cite{PL}. We
shall prove a more precise result.

\begin{Th}\label{diamUB} Let $E$ be a Banach space and
 let $A\subset E$  a closed and uniformly convex subset with modulus
$\delta_{A} (\cdot)$. Then for any $\ep \in (0,\diam A)$
$$
\diam A\le \left(\left[
\frac{\ep}{\delta_{A}(\ep)}\right]+1\right)\cdot\ep,
$$
where $[x]$ is the largest integer $\le$ $x$.
\end{Th}

\proof For any unit vector $p\in b(A)$ and
any $t>0$ we
define a convex closed set:
$$
A_{p}(t) = A\cap \{ x\in E\ |\ \l p,x\r\ge s(p,A)-t \}.
$$
We obtain
from the definition of the supporting function  that
$A_{p}(t)\ne\emptyset$ for any $t>0$, $p\in b(A)$, $\| p \|=1$,
and if $0<t_{1}<t_{2}$ then $A_{p}(t_{1})\subset A_{p}(t_{2})$.

We shall show that for any unit $p\in b(A)$ the following holds
$$
\lim\limits_{t\to+0}\diam A_{p}(t)=0.
$$
Suppose that for some unit $p\in b(A)$ there exist $d>0$ and
$t_{k}\to +0$ with $\diam A_{p}(t_{k})\ge d$. The latter means
that there exist points $x_{k}^{1}$, $x_{k}^{2}$ from
$A_{p}(t_{k})$ with $\| x_{k}^{1}-x_{k}^{2}\|>d/2,$ for all $k$.
It follows from uniform convexity of the set $A$ that
$$
B_{\delta_{A}\left(\frac{d}{2}\right)}\left(
\frac{x_{k}^{1}+x_{k}^{2}}{2}\right)\subset A.
$$
However, by taking the supporting functions of the sets from this
inclusion we obtain the following:
$$
s\left( p,B_{\delta_{A}\left(\frac{d}{2}\right)}\left(
\frac{x_{k}^{1}+x_{k}^{2}}{2}\right)\right) = \frac12 (\l
p,x_{k}^{1}\r +\l p,x_{k}^{2}\r) +
\delta_{A}\left(\frac{d}{2}\right)\ge s(p,A)-t_{k}+
\delta_{A}\left(\frac{d}{2}\right)\ge s(p,A).
$$
The last inequality holds for sufficiently large $k$ (when
$\delta_{A}\left(\frac{d}{2}\right)>t_{k}$). This contradiction
shows that $\diam A_{p}(t)\to 0$, $t\to+0$.

By the completness of $A$
we conclude that
$$
\bigcap\limits_{t>0}A_{p}(t)=\{ a(p)\}.
$$

We have thus
proved that for any unit vector $p\in b(A)$ there
exists $a(p)=\arg\max\limits_{x\in A}\l p,x\r$.

\def\cone {\mbox{cone}\,}

Let's fix arbitrary points $x,y\in\d A$. By the separation theorem
there exist unit vectors $q_{1},q_{2}\in E^{*}$ such that $\l
q_{1},x\r=s(q_{1},A)$, $\l q_{2},y\r=s(q_{2},A)$. If $q_{1}\ne
-q_{2}$ then let $D=\d B_{1}^{*}(0)\cap\cone\{ q_{1},q_{2}\}$. If
$q_{1}=-q_{2}$ then choose any $q_{3}\in \d B_{1}^{*}(0)$ with
$s(q_{3},A)<+\infty$ and define $D=\d B_{1}^{*}(0)\cap\cone\{
q_{1},q_{2},q_{3}\}$. Note that for any $q\in D$,
$s(q,A)<+\infty$.

By \cite[Theorem 11.9]{Leht} for any 2-dimen\-si\-onal subspace
$\L\subset E^{*}$ the length of the curve $\L\cap\d B_{1}^{*}(0)$
is less than $8$ (in the $\|\cdot \|_{*}$-norm). Thus the length
of $D$ is less than $4$.
Choose $N=\left[ \frac{\ep}{\delta_{A}(\ep)}\right]+1$ and points
$\{ p_{i}\}_{i=0}^{N}\in D$ which decompose
the length of $D$ into
$N$ equal parts; $p_{0}=q_{1}$, $p_{N}=q_{2}$ and $\|
p_{i-1}-p_{i}\| <\frac{4}{N}\le \varphi (\ep)$, $i=1,\ldots,N$.

By the previous considerations  we obtain that for all $1\le i\le
N-1$ there exists  $x_{i}\in \d A$ with $\l
p_{i},x_{i}\r=s(p_{i}, A)$.
By  Lemma 2.\ref{2} we have $\|  x_{i-1}- x_{i}\|<\ep$ and
$$
\| x-y\|\le \sum\limits_{i=1}^{N}\|  x_{i-1}- x_{i}\|\le \ep\cdot
N.
$$
The points $x,y$ are arbitrary boundary points of $A$, hence
$\diam A\le \ep\cdot N$.\qed

\begin{Corol}\label{C2}
By Theorem 2.\ref{diamUB} we have
$$
\delta_{A}(\ep )\le\frac{\ep^{2}}{\diam A-\ep},\quad\forall
\ep\in (0,\diam A).
$$
This means that $\delta_{A}(\ep)\le C\cdot \ep^{2}$ for any
convex closed bonded set $A$.
\end{Corol}

For balls this statement follows from the well-known
Day-Nordlendar theorem \cite{Diestel} which asserts
that if $E$ is
a Banach space then the modulus of convexity for $E$, i.e. the
modulus of convexity for the unit ball, satisfies the estimate
$\delta_{E}(\ep)\le 1-\sqrt{1-\frac{\ep^{2}}{4}}$, $\forall\ep\in
(0,2)$.

Next we shall prove a result which is very close to the
Day-Nordlendar theorem.
%-----------------------------------------------------------------------------

\begin{Th}\label{Day+Nordlendar}
Let $E$ be a Banach space and $A\subset E$  a closed and
uniformly convex set with modulus $\delta_{A}(\cdot)$, $\diam
A=1$. Let $r_{0}>0$ and $a\in E$ be such that $B_{r_{0}}(a)\subset
A$. Then for all $\ep\in (0,1) $:
\begin{equation}\label{DN1}
\delta_{A}(2 r_{0}\ep)\le \frac12\left(1-\sqrt{1-\ep^{2}}\right).
\end{equation}

 In
(2.\ref{DN1})
the equality
takes place when $A$ is the Euclidean ball of
diameter 1 in the Euclidean space (with $r_{0}=\frac12$).
\end{Th}

\proof Without loss of generality we can assume that $a=0$.
Let $B=A\cap (-A)$. Note that the set $B$ is bounded, has a
nonempty interior ($B_{r_{0}}(0)\subset B$) and its center of
symmetry in zero. Hence we can consider the set  $B$ as the ball
of radius $\frac12$ and we have:
\begin{equation}\label{DN2}
 B_{r_{0}}(0)\subset B\subset
B_{\frac12}(0).
\end{equation}
Let's say few words about the second inclusion in (2.\ref{DN2}).
If $x\in B$, then $-x\in B$, and $2\|x\|=\|x-(-x)\|\le \diam B=1$.
Therefore $B\subset B_{\frac12}(0)$.
  By $\|\cdot\|_{B}$ we denote the
new norm with the unit ball $2B$.

For any convex closed bounded set  $C\subset E$ we shall consider
the modulus of convexity:
$$
\delta_{C}^{B}(\ep)= \sup\left\{ \delta\ge 0\ \left|\ \delta\cdot
2 B+ \frac{x_{1}+x_{2}}{2}\right.\subset C,\ \forall
x_{1},x_{2}\in \d C:\ \| x_{1}-x_{2}\|_{B}=\ep \right\}.
$$

Let $x_{1},x_{2}\in B$ and $\| x_{1}-x_{2}\|_{B}=\ep\in (0,1)$.
From $\delta_{A}^{B}=\delta_{-A}^{B}$ we have
$$
\frac{x_{1}+x_{2}}{2}+ 2B\delta_{A}^{B}(\ep)\subset A,\qquad
\frac{x_{1}+x_{2}}{2}+ 2B\delta_{-A}^{B}(\ep) =
\frac{x_{1}+x_{2}}{2}+ 2B\delta_{A}^{B}(\ep)\subset -A.
$$
By definition, $B=A\cap (-A)$, so we obtain that
$\frac{x_{1}+x_{2}}{2}+\delta_{A}^{B}(\ep)\cdot 2B\subset B$ and
thus $\delta_{A}^{B}(\ep)\le \delta_{B}^{B}(\ep)$ for all $\ep \in
(0,1)$. From the equality
$\delta_{B}^{B}(\ep)=\frac12\delta_{2B}^{B}(2\ep)$, using
Day-Nordlendar theorem \cite[Theorem  3.3.1]{Diestel} for the unit
ball $2B$, we obtain for all $\ep\in (0,1)$
$$
\delta_{B}^{B}(\ep)=\frac12\delta_{2B}^{B}(2\ep)\le
\frac12\left(1-\sqrt{1-\frac{(2\ep)^{2}}{4}}\right),
$$
and $\delta_{A}^{B}(\ep)\le
\frac12\left(1-\sqrt{1-\ep^{2}}\right)$ for all $\ep\in (0,1)$.

We conclude from inclusions (2.\ref{DN2}), that for any
$x_{1},x_{2}\in E$ the inequalities $2 r_{0}\|
x_{1}-x_{2}\|_{B}\le \| x_{1}-x_{2}\|\le \| x_{1}-x_{2}\|_{B}$
hold. If $x_{1},x_{2}\in \d A$, $\| x_{1}-x_{2}\|_{B}=\ep$ and
$\ep\in (0,1)$, then $\delta_{A}(2 r_{0} \ep )\le \delta_{A}(\|
x_{1}-x_{2}\|)$. Since for any $\delta\ge 0$ the condition
$\frac{x_{1}+x_{2}}{2}+\delta B_{1}(0)\subset A$ implies the
condition $\frac{x_{1}+x_{2}}{2}+\delta\cdot 2B\subset A$, it
follows that $\delta_{A}(\| x_{1}-x_{2}\| )\le
\delta_{A}^{B}(\ep)$. Therefore we get the formula (2.\ref{DN1}).

An easy calculation show that in the case when the set  $A$ is a
Euclidean ball of diameter 1 in the
Euclidean space with
$r_{0}=\frac12$ we get the equality in the formula (2.\ref{DN1}).
 \qed

%-------------------------------------------------------------------------------------------DN

\begin{Th}\label{Uninorm}
In every Banach space $E$ there exists a closed uniformly convex
set $A$ if and only if the space $E$ admits an equivalent
uniformly convex norm.
\end{Th}

\proof Due to Theorem 2.\ref{diamUB} we must consider only bounded
sets.
If the space $E$ admits an equivalent uniformly convex norm then
the unit ball of this norm is a uniformly convex set. Let us
prove the converse statement.

Let $A\subset E$ be closed and uniformly convex set with modulus
$\delta_{A}$. Suppose that $0\in \i A$. As we can see from the
proof of Theorem 2.\ref{Day+Nordlendar}, the set $B=A\cap (-A)$
is a uniformly convex ball of equivalent norm.\qed

Note that a Banach space which is equivalent to a uniformly convex
space, is reflexive \cite{Diestel}. Thus we can further use
reflexivity without loss of generality. The reflexivity of the
Banach space with bounded nonsingleton uniformly convex set was
mentioned in \cite{Polyak}.
We also note that nonreflexive spaces (e.g., the spaces
$C([0,1])$, $L_{1}([0,1])$, $L_{\infty}([0,1])$, $l_{1}$,
$l_{\infty}$) do not contain uniformly convex sets.

Recall that in any finite-dimensional Banach space the class of strictly convex compacta
coincides with the class of uniformly convex sets. This fact
easily follows from compactness of sets from two classes. It is
well-known \cite{Diestel} that in infinite-dimension spaces there
exist strictly but nonuniformly convex balls.

We wish to mention an
important class of
uniformly convex sets.
Let $E$ be a uniformly convex Banach space. The set $A\subset E$
is strongly convex with radius $R>0$ \cite[Chapters 3,
4]{Polovinkin&Balashov} (or $R$-convex \cite{Olech+Fr}) if
$A=\bigcap\limits_{x\in X}B_{R}(x)\ne\emptyset$, where $X\subset
E$ an arbitrary subset. It is easy to see that the modulus of
convexity for $A$ is $\delta_{A}(\ep)\ge
R\delta_{E}\left(\frac{\ep}{R}\right)$ for all $\ep\in (0,\diam
A)$. Here $\delta_{E}$ is the
modulus of convexity for the space $E$.

\section{Applications to the set-valued analysis and \\ the splitting problem for selections}

Let $\{F(t)\}_{t\in T}$ be any collection of convex closed sets
and let $\diam F(t)\ge r_{0}>0$ for all $t$. Suppose that each
set $F(t)$ is uniformly convex with modulus $\delta_{t}(\ep)$.
Then under the assumption that $\delta(\ep )=\inf\limits_{t\in
T}\delta_{t}(\ep)>0$ for all $\ep\in (0,r_{0})$, the set
$F=\cap_{t\in T}F(t)$ is uniformly convex with modulus
$\delta_{F}(\ep)\ge \delta (\ep)$ for all $\ep\in (0,r_{0})$
(this set can
also be empty or a singleton). Note that Lemmata
2.1, 2.2 and Theorem 2.1 are valid for the function $\delta (\ep)$
and the set $F$.

Consider as an example the set $A$, which can be represented as
the intersection of closed balls of radius 1 in Hilbert space
$\H$. The modulus of convexity for the unit ball from $\H$ is
$\delta_{\H}(\ep)=1-\sqrt{1-\frac{\ep^{2}}{4}}\ge
\frac{\ep^{2}}{8}$ for all $\ep\in (0,2)$ and $\delta_{A}(\ep)\ge
\delta_{\H}(\ep)$. By Corollary 2.2 we have that $\varphi
(\ep)\ge\ep/2$ and $\| x_{1}-x_{2}\|\le 2\| p_{1}-p_{2} \|$. So
we conclude that the gradient $\nabla s(p,A) =
\arg\max\limits_{x\in A}\l p,x\r$ of supporting function for the
set $A$ is a Lipschitz function with respect to $p$. This result
was proved in \cite{Polovinkin&Balashov} by different methods.

Next we shall consider set-valued mappings $F:T\to
2^{E}\backslash\emptyset$ from a metric space $(T,\rho)$ to a
Banach space $E$. Suppose that there exists $r_{0}>0$ such that
for any $t\in T$ we can find a point $a(t)\in E$ with
$B_{r_{0}}(a(t))\subset F(t)$. Suppose that any set $F(t)$ is
closed and uniformly convex with modulus $\delta_{t}(\ep)$,
$\ep\in (0,\diam F(t))$. If $\delta (\ep)=\inf\limits_{t\in
T}\delta_{t}(\ep)>0$ for all $\ep\in (0,2r_{0}]$ then we say that
the images $F(t)$, $t\in T$, are uniformly convex with modulus
$\delta (\ep)$, $\ep\in (0,2r_{0}]$. It's easy to see that Lemmata
2.1, 2.2 and Theorem 2.1 are valid for any set $F(t)$ when instead
of the modulus $\delta_{F(t)}$ we take the modulus $\delta$.

For an increasing function $\delta : [0,d]\to [0,\Delta]$ we
define the inverse function $\delta^{-1}$ as follows: for
$x_{0}\in [0,\Delta]$ let $\delta^{-1}(x_{0})=y_{0}\in [0,d]$.
Here $\delta(y_{0}-0)\le x_{0}\le \delta (y_{0}+0)$; $\delta
(y_{0}\pm 0)=\lim\limits_{y\to y_{0}\pm 0}\delta (y)$. Note that
the function $\delta^{-1}$ is continuous on the segment
$[0,\Delta]$.

We shall use $\conv A$ to denote the convex hull of the set $A$.
The Hausdorff distance $h(A,B)$ between
sets $A$ and $B$ in a
Banach space $E$ is defined as follows:
$$
h(A,B)=\inf\{ r>0\ |\ A\subset B+B_{r}(0),\ B\subset A+B_{r}(0)\}.
$$

\begin{Th}\label{22}
Let $(T,\rho)$ be a metric space and $E$ a reflexive Banach
space. Suppose that the set-valued mappings $F_{i}:T\to
2^{E}\backslash\emptyset$, $i=1,2$, have convex closed images.
Let $F_{i}$, $i=1,2$, be uniformly continuous in the Hausdorff
metric, i.e. there exist nonnegative infinitely small at zero
numerical functions $\omega_{i}$, such that for all
$t_{1},t_{2}\in T$ we have the following:
$$
h(F_{i}(t_{1}),F_{i}(t_{2}))\le \omega _{i}(\rho (t_{1},t_{2})).
$$

Let the images $F_{1}(t)$ be uniformly convex with modulus $\delta
(\ep)$, $\ep\in (0,2r_{0}]$. Let $\Delta_{0}=\delta (2r_{0})$ and
 $H(t)=F_{1}(t)\cap F_{2}(t)\ne\emptyset$ for all $t\in T$.

Then
\begin{equation}\label{f5}
 h(H(t_{1}),H(t_{2}))\le \omega_{1}(\rho (t_{1},t_{2})) +
2\omega_{2}(\rho (t_{1},t_{2})) + f\Bigl(\omega_{1}(\rho
(t_{1},t_{2})) + \omega_{2}(\rho (t_{1},t_{2}))\Bigr),
\end{equation}
 where
\begin{equation}\label{f6}
f(x)=\left\{
\begin{array}{c}
\delta^{-1}\left(\frac{x}{2}\right),\qquad x<2\Delta_{0},\\
\frac{Mx}{2\Delta_{0}},\qquad x\ge 2\Delta_{0}
\end{array}\right. ,
\end{equation}
 and $M =\sup\limits_{t\in T}\diam F_{1}(t)\le r_{0}\left( \left[
\frac{r_{0}}{\delta (r_{0})} \right]+1\right)$.
\end{Th}

\proof We define $\omega_{1}=\omega_{1}(\rho (t_{1},t_{2}))$,
$\omega_{2}=\omega_{2}(\rho (t_{1},t_{2}))$.
 Let $b_{1}\in H(t_{1})$. Let's fix $k>1$. We shall prove
that there exists point $a(t_{2})\in H(t_{2})$ such that

\begin{equation}\label{ffn}
 \| a(t_{2})-b_{1}\|\le
f(\omega_{1}+k\omega_{2})+\omega_{1}+2k\omega_{2}.
\end{equation}

We obtain
from formula (3.\ref{ffn})  the following:
$$
h(H(t_{1}),H(t_{2}))\le
f(\omega_{1}+k\omega_{2})+\omega_{1}+2k\omega_{2},
$$
and keeping in mind that the function $f$ is continuous from the
right (see (3.\ref{f6})), we  take the limit $k\to 1+0$ and
obtain formula (3.\ref{f5}).

Let $b(t_{2})\in F_{2}(t_{2})$: $\| b(t_{2})-b_{1}\|\le
kh(F_{2}(t_{2}),F_{2}(t_{1}))\le k\omega_{2}$.
If $b(t_{2})\in F_{1}(t_{2})$ then we can take $a(t_{2})=b(t_{2})$
and we conclude that formula (3.\ref{ffn}) is valid. Further we
shall assume
that $b(t_{2})\notin F_{1}(t_{2})$.

Let $c(t_{2})\in H(t_{2})\subset F_{1}(t_{2})$. Let
$b_{\pi}(t_{2})$ be the metric projection of the point $b(t_{2})$
onto $F_{1}(t_{2})$. The point $b_{\pi}(t_{2})$ exists because
the space  $E$ is reflexive. Consider
the point $a(t_{2})$
which is the nearest to the point $b(t_{2})$ of the set
$F_{1}(t_{2})\cap\conv\left\{ b(t_{2}),c(t_{2}) \right\}$. By
definition, $a(t_{2})\in F_{1}(t_{2})$ and $a(t_{2})\in
\conv\left\{ b(t_{2}),c(t_{2}) \right\}\subset F_{2}(t_{2})$.
This implies that $a(t_{2})\in H(t_{2})$.

Let $z(t_{2})=\frac{a(t_{2})+b(t_{2})}{2}$, $\tilde z (t_{2})=
\frac{a(t_{2})+b_{\pi}(t_{2})}{2}$. Since
$$
\| z(t_{2}) - \tilde z (t_{2})\|=\frac12 \|
b(t_{2})-b_{\pi}(t_{2})\|,\qquad  B_{\delta (\|
a(t_{2})-b_{\pi}(t_{2})\|)}(\tilde z (t_{2}))\subset F_{1}(t_{2}),
$$
it follows from the condition $z(t_{2})\notin F_{1}(t_{2})$
that
\begin{equation}\label{f7}
\delta (\| a(t_{2})-b_{\pi}(t_{2})\|)\le  \| z(t_{2}) - \tilde z
(t_{2})\| = \frac12 \| b(t_{2})-b_{\pi}(t_{2})\|.
\end{equation}
 So we have following estimate:
$$
\begin{array}{l}
\| b(t_{2})-b_{\pi}(t_{2})\| = \rho (b(t_{2}),F_{1}(t_{2}))\le
\rho (b_{1},F_{1}(t_{2}))+\| b(t_{2})-b_{1}\| \le\\
\qquad\qquad\le h(F_{1}(t_{1}),F_{1}(t_{2}))+k\omega_{2}\le
\omega_{1}+k\omega_{2}.
\end{array}
$$
By the last formula and by (3.\ref{f7}) we have that $\delta (\|
a(t_{2})-b_{\pi}(t_{2})\|)\le\frac12 (\omega_{1}+k\omega_{2})$.

If $\omega_{1}+k\omega_{2}<2\Delta_{0}$ then
$$
\| a(t_{2})-b_{\pi}(t_{2})\| \le\delta^{-1}\left( \frac12
(\omega_{1}+k\omega_{2})\right).
$$
If $\omega_{1}+k\omega_{2}\ge 2\Delta_{0}$ then
$$
\| a(t_{2})-b_{\pi}(t_{2})\| \le
\frac{\omega_{1}+k\omega_{2}}{2\Delta_{0}}M.
$$
Thus in both cases we have $\| a(t_{2})-b_{\pi}(t_{2})\|\le
f(\omega_{1}+k\omega_{2})$.
 Finally,
$$
\| a(t_{2})-b_{1}\| \le \| a(t_{2})-b_{\pi}(t_{2})\|+\|
b_{\pi}(t_{2})-b(t_{2})\|+\| b(t_{2})-b_{1}\|\le
f(\omega_{1}+k\omega_{2})+\omega_{1}+2k\omega_{2}.
$$
\qed

Theorem 3.\ref{22} has important consequences. It follows from
Corollary 2.3 that the modulus of convexity  $\delta (\ep)$ of
sets $F_{1}(t)$ in Theorem 3.\ref{22} does not exceed  $C\cdot
\ep^{2}$. Hence the H\"older condition with the power no greater
than
$\frac12$ with respect to the Hausdorff metric is typical for the
product of intersections of two Lipschits set-valued mappings. We
need to invoke good mutual geometric properties of $F_{1}$ and
$F_{2}$ if we want to obtain power greater
than $\frac12$ (see for
example \cite[Theorem 2.2.1]{Polovinkin&Balashov}). Under the
conditions of Theorem 3.\ref{22} the result is the best possible.

\Example \rm In the Euclidean plane $\R^{2}$ with the standard
basis $x_{1}Ox_{2}$ we consider (for $t\ge 0$)
$$
F_{1}(t)=F_{1}=\{ (x_{1},x_{2})\ |\  x_{2}\ge |x_{1}|^{p}\}\bigcap
B_{1}(0),\ p\ge 2,\qquad F_{2}(t)=\{ (x_{1},t)\ |\ x_{1}\in\R\}.
$$
It is easy to see that if  $\ep>0$ is sufficiently small then the
modulus of convexity $F_{1}$ equals $\delta_{1}(\ep
)=\frac{\ep^{p}}{2^{p}}$ (and it realized on the segment
$\left[\left( -\frac{\ep}{2},\frac{\ep^{p}}{2^{p}}\right),\left(
\frac{\ep}{2},\frac{\ep^{p}}{2^{p}} \right) \right]$). The
intersection of $F_{1}(t)$ and $F_{2}(t)$ is
$H(t)=[-t^{1/p},t^{1/p}]\times\{ t\}$. Let $t_{1}>0$,
$t_{2}=2t_{1}$. Then
$$
h(H(t_{1}),H(t_{2}))\ge
(2t_{1})^{1/p}-t_{1}^{1/p}=(2^{1/p}-1)\cdot |t_{2}-t_{1}|^{1/p}=
\frac{2^{1/p}-1}{2}\,\delta_{1}^{-1}(|t_{2}-t_{1}|).
$$

\Example \rm Consider the following extremal problem
\begin{equation}\label{f8}
 \min\limits_{x\in A}g(x).
\end{equation}
Suppose that the function $g$ has closed and uniformly convex
level sets ${\mathcal L}_{g}(\beta)=\{ x\in E\ |\ g(x)\le \beta
\}$. The function $g$ itself
cannot be  convex. We shall consider
two problems (3.\ref{f8}) with the same function and convex
closed sets $A_{i}$, $i=1,2$. Suppose that the point $u_{i}$ is
the solution of the problem (3.\ref{f8}) with the set $A=A_{i}$,
i.e. $\{ u_{i}\}=A_{i}\cap{\mathcal L}_{g}\left( \min\limits_{x\in
A_{i}}g(x)\right)$. We shall estimate the value $\| u_{1}-u_{2}\|$
through the distance $h=h(A_{1},A_{2})$.

Note that for convex functions and sets such problems were
considered e.g., in \cite{Berd}, \cite{Mar}.
Let $g(u_{1})\le g(u_{2})$. Then
$$
u_{1}=A_{1}\cap{\mathcal L}_{g}(g(u_{1}))\subset
A_{1}\cap{\mathcal L}_{g}(g(u_{2})).
$$
Let the set ${\mathcal L}_{g}(g(u_{2}))$ be uniformly convex with
modulus $\delta$. Let $F_{1}(A)={\mathcal L}_{g}(g(u_{2}))$ be a
constant mapping with the modulus of continuity $\omega_{1}=0$,
and let $F_{2}(A)=A$ be a mapping with the modulus of continuity
$\omega_{2}(t)=t$. By Theorem 3.\ref{22} we have
$$
h\left(F_{1}(A_{1})\cap F_{2}(A_{1}), F_{1}(A_{2})\cap
F_{2}(A_{2})\right)\le 2h+f(h),
$$
where function  $f$ is defined in (3.\ref{f6}) and $M=\diam
{\mathcal L}_{g}(g(u_{2}))$, $\Delta_{0}=\lim\limits_{\ep\to
\diam {\mathcal L}_{g}(g(u_{2}))-0}\delta (\ep )$.
Therefore for all
$t>1$
$$
A_{1}\cap{\mathcal L}_{g}(g(u_{2}))\subset A_{2}\cap{\mathcal
L}_{g}(g(u_{2}))+t(2h+f(h))B_{1}(0)=u_{2}+t(2h+f(h))B_{1}(0),
$$
i.e.
\begin{equation}\label{f9}
\| u_{1}-u_{2}\|\le 2h(A_{1},A_{2})+f(h(A_{1},A_{2})).
\end{equation}

We now consider applications of the above results to the splitting
problem for selections.

\Example \rm (Question 4.6 from \cite{Repovs+Semenov}). Do there
exist for every closed convex sets $A$, $B$ and $C=A+B$
continuous functions $a:C\to A$ and $b:C\to B$ with the property
that $a(c)+b(c)=c$ for all $c\in C$? Similar questions were also
considered in previous papers, see \cite{LM-LS} for details.

\begin{Lm}
Let the space  $E$ be uniformly convex with modulus $\delta_{E}$.
Let $A\subset E$ be a closed and uniformly convex set with modulus
$\delta_{A}$, and $B\subset E$ a convex and closed set. Then there
exist uniformly continuous functions $a:C\to A$ and $b:C\to B$
such that $a(c)+b(c)=c,$ for all points $c\in C$.
\end{Lm}

\proof Suppose that $0\notin A$. For any $c\in C$ we define sets
$F_{1}(c)=A$, $F_{2}(c)=c-B$. Then (in terms of Theorem
3.\ref{22}) $\omega_{1}=0$, $\omega_{2}(t)=t$, $F_{1}$ has
uniformly convex images with modulus $\delta_{A}$. Note that
$H(c)=(c-B)\cap A$ is nonempty for all $c\in C$.

Let's define
$M=\diam A$, $\Delta_{0}=\lim\limits_{\ep\to\diam A-0
}\delta_{A}(\ep)$. By Theorem 3.\ref{22}
\begin{equation}\label{f12}
 h(H(c_{1}),H(c_{2}))\le 2\| c_{1}-c_{2}\|+f(\|
c_{1}-c_{2}\|),
\end{equation}
 where $f$ is from (3.\ref{f6}).

Let $r=\inf\limits_{a\in A}\| a\|>0$, $R=\sup\limits_{a\in A}\|
a\|$. All balls $B_{t}(0)$, $t\in [r,R]$, are uniformly convex
with modulus $\delta(\ep)=R\delta_{E}(\frac{\ep}{R})$, $\ep \in
(0,2r]$. Let $a(c)=\arg\min\limits_{x\in H(c)}\| x\|$. Let's
define $\Delta_{E}=\delta(2r)$,
$$
f_{E}(t)=\left\{
\begin{array}{c}
\delta^{-1}\left(\frac{t}{2}\right),\qquad t<2\Delta_{E},\\
\frac{Rt}{\Delta_{E}},\qquad t\ge 2\Delta_{E}.
\end{array}\right.%\eqno(11)
$$
Using (3.\ref{f9}) from Example 3.2 and (3.\ref{f12}) we have
$$
\begin{array}{l}
\| a(c_{1})-a(c_{2})\|\le
2h(H(c_{1}),H(c_{2}))+f_{E}(h(H(c_{1}),H(c_{2})))\le \\ \qquad\le
4\| c_{1}-c_{2}\|+2f(\| c_{1}-c_{2}\|)+f_{E}\left( 2\|
c_{1}-c_{2}\|+f(\| c_{1}-c_{2}\|)\right).
\end{array}
$$
So we have built uniformly continuous selections $a(c)\in
H(c)\subset A$ and $b(c)=c-a(c)\in B$.\qed

\begin{Remark}\label{InRn}
Note that in the case $E=\R^{n}$ we can define $a(c)$ as
$a(c)=s(H(c))$, where $s(H(c))$ is the Steiner point of the set
$H(c)$. The Steiner point is a Lipschitz selection of convex
compacta from $\R^{n}$ with the Lipschitz constant
$L_{n}=\frac{2}{\sqrt{\pi }}\frac{\Gamma\left(
\frac{n}{2}+1\right)}{\Gamma\left( \frac{n+1}{2}\right)}$
\cite{Aubin+Frankowska,Polovinkin&Balashov}. From this and
by
formula (3.\ref{f12}) we get
$$
\| a(c_{1})-a(c_{2})\| \le L_{n}\cdot\left( 2\| c_{1}-c_{2}\|+f(\|
c_{1}-c_{2}\|)\right).
$$
\end{Remark}

\begin{Remark}\label{InUniE}
Let $A$ and $B$ be closed convex subsets of the reflexive Banach
space $E$ and let the set $A$ be strictly convex and bound. Let
$0\in \Int A$.

 Let $c\in A+B$, $\varrho_{A} (c,B)=\inf\{ t>0\ |\ c\in B+tA\}$,
 and
$$
b(c)=\left( c-\varrho_{A}(c,B)A\right)\cap B.
$$
The set $b(c)$ is a point.
This follows from the reflexivity of the
space $E$ (the set $B+tA$ is closed for all $t\ge 0$) and
strictly convexity of the set $A$. The point $b(c)$ is projection
of the point $c$ in the sense of the set $A$ on the set $B$. Note
that in above situation $\varrho_{A}(c,B)\in [0,1]$.

If this projection $b(c)$ uniformly continuously depends on $c$,
then $b(c)\in B$
is a
uniformly continuous selection of $B$ and $a(c)=c-b(c)\in \varrho_{A}(c,B)A\subset
A$
is a uniformly continuous selection of $A$.

In particular, if the spaces $E$ and $E^{*} $ have moduli of
convexity of the second order and $A=B_{1}(0)$ then by the
results from
\cite{Alber} we obtain
that the
projection $b(c)$ satisfies the
Lipschitz
condition with respect to
$c$. In particular, this takes place in
the Hilbert space. It would be very interesting to describe all
spaces and pairs of sets ($A$ and $B$) for which the projection
$b(c)$ of the point $c$ in the sense of the set $A$ on the set $B$
satisfies the
Lipschitz condition.

\end{Remark}

\Example

Hereafter, the sum of Banach spaces $E_{1}\oplus E_{2}$ will be
defined as follows: $w=(u,v)\in E_{1}\oplus E_{2}$, $\|
w\|=\max\{ \| u \|_{E_{1}}, \| v \|_{E_{2}}\}$.

\begin{Lm}\label{sec}
Let $T$ be a metric space, $E_{i}$ a reflexive Banach spaces, and
$F_{i}:T\to 2^{E_{i}}$ uniformly continuous set-valued mappings
with modulus of continuity $\omega$, i.e.
$$
h((F_{1}(t_{1}),F_{2}(t_{2})),(F_{1}(t_{1}),F_{2}(t_{2})))\le
\omega (\rho (t_{1},t_{2})),\qquad \forall t_{1},t_{2}\in T, \
i=1,2.
$$
Suppose that the images
$F_{i}(t)$ are uniformly convex sets with modulus
$\delta (\ep)$, $i=1,2$ and $\ep\in (0,2r_{0}]$;
$\Delta_{0}=\delta(2r_{0})$. Let $\L\subset E_{1}\oplus E_{2}$ be
a closed subspace and suppose that
there exists $C>0$
such that for any
$w_{1}=(u_{1},v_{1})\in \L$, $w_{2}=(u_{2},v_{2})\in \L$ we have
$\| u_{1}-u_{2}\|_{E_{1}}\ge C\| w_{1}-w_{2}\|$ and  $\|
v_{1}-v_{2}\|_{E_{2}}\ge C\| w_{1}-w_{2}\|$ (i.e. $\L$ is not
"parallel" to $E_{1}$ and $E_{2}$).

Let $M=\sup\limits_{t\in T}\diam (F_{1}(t),F_{2}(t))<+\infty$.
Define the set-valued map $H(t)=(F_{1}(t),F_{2}(t))\cap\L\ne\emptyset$
for all $t\in T$.
Then
\begin{equation}\label{USP}
 h(H(t_{1}),H(t_{2}))\le
 \omega (\rho
(t_{1},t_{2}))+ \frac{1}{C}f(\omega (\rho
(t_{1},t_{2}))),\qquad\forall t_{1},t_{2}\in T,
\end{equation}
where the function $f$ is from formula (3.\ref{f6}).
\end{Lm}

\proof Let $w_{0}\in H(t_{0})$. Let's fix $k>1$. We shall prove
that there exists a point a $w_{1}\in H(t)$ with the following property:
$$
\| w_{0}-w_{1}\|\le k\omega (\rho (t_{0},t))+ \frac{1}{C}f(k\omega
(\rho (t_{0},t))).
$$
Thus
$$h(H(t_{1}),H(t_{2}))\le k\omega (\rho (t_{0},t))+
\frac{1}{C}f(k\omega (\rho (t_{0},t)))$$
and we obtain
(3.\ref{USP})
by
taking the limit
$k\to 1+0$.

Let $w\in (F_{1}(t),F_{2}(t))$ be a point
such that $\|
w_{0}-w\|\le k\omega (\rho (t_{0},t))$.
%Let $\tilde w$ be the
%metric projection $w$ onto $\L$ (exists from theorem
%\ref{Uninorm}).
%Let $B=B_{\|\tilde w-w_{0} \|}(w_{0})\cap\L$.
Define $w_{1}\in H(t)$ to be the point from the set $H(t)$ which
is the nearest to the point $w_{0}$ ($w_{1}$ exists by the
reflexivity of $E_{i}$, $i=1,2$).

%If $w_{1}\in B$ then
%$$
%\| w_{0}-w_{1}\|\le \| w_{0}-w\|+\|w-w_{1} \|\le \| w_{0}-w\|+\|
%w-\tilde w\|\le 2\| w_{0}-w\| = 2\omega (\rho (t_{0},t))
%$$
%and inequality (\ref{USP}) holds true.

%If $w_{1}\notin B$ then
Let $w_{2}=\frac12 (w+w_{1})$.  If $z\in \L$ is the middle point
of the segment $[w_{1},w_{0}]$ then
$$
\| w_{2}-z\| = \frac12 \| w-w_{0}\|\le k\frac12 \omega (\rho
(t_{0},t)).
$$
Thus we must require
$\delta (C\| w-w_{1}\|)\le k\frac12 \omega
(\rho (t_{0},t))$. Otherwise we would have, since $\L$
is  "parallel"
neither to $E_{1}$
nor to $E_{2}$,
the following
contradiction:
$$
z\in B^{E_{1}\oplus E_{2}}_{\delta (C\| w-w_{1}\|)}(w_{2})\cap
\L\subset (F_{1}(t),F_{2}(t))\cap\L = H(t),
$$
with the inequality $\| w_{1}-w_{0}\|\le\| z-w_{0}\|$.

If $k\omega (\rho (t_{0},t))<2\Delta_{0}$ then
$$\| w-w_{1}\|\le
\frac{1}{C}\delta^{-1}\left( k\frac12\omega (\rho
(t_{0},t))\right).$$
If $k\omega (\rho (t_{0},t))\ge 2\Delta_{0}$
then
$$\| w-w_{1}\|\le \frac{1}{C}\frac{k\omega (\rho
(t_{0},t))}{2\Delta_{0}}M.$$

In both cases $\| w-w_{1}\|\le \frac{1}{C}f(k\omega (\rho
(t_{0},t)))$.
Finally,
$$
\| w_{0}-w_{1}\|\le \| w_{0}-w\|+\| w-w_{1}\|\le k\omega (\rho
(t_{0},t)) + \frac{1}{C}f(k\omega (\rho (t_{0},t))).
$$
\qed

\begin{Remark}\label{best}
The result (3.\ref{USP}) of Lemma 3.\ref{sec} is exact. Let $T$
be the space of convex closed bounded subsets of the Hilbert space
$\H$ with the Hausdorff distance, $E=\H$.
Define set-valued mappings $F_{i}:T\to 2^{\H}$, $i=1,2$, as
follows:
$$
\forall A\in T\qquad F_{1}(A)=A,\quad
F_{2}(A)=B_{\varrho(0,A)}(0),
$$
where $\varrho(0,A)=\inf\limits_{a\in A}\| a\|$. Note that
$\delta_{F_{2}}(\ep)=C\cdot\ep^{2}$ (the modulus of convexity for the
Hilbert space). Obviously, $T\ni A\to F_{i}(A)$, $i=1,2$, are
Lipschitz functions in the Hausdorff metric.

Let $L:\H\oplus \H\to \H$, $L(y_{1},y_{2})=y_{1}-y_{2}$, $\L=\ker
L =  \{ (y_{1},y_{2})\in\H\oplus \H\ |\ y_{1}-y_{2}=0\}$. Then
$$
(F_{1}(A),F_{2}(A))\cap\L = \{ (p(A),p(A))\},
$$
where $p(A)$ is the metric projection of the zero on the set $A$.
It follows by well-known results of Daniel \cite{Daniel}, that
$T\ni A\to p(A)$ is a H\"older function with power $\frac12$ in
the Hausdorff metric.

\end{Remark}

\begin{Th}
Let $T$ be a metric space, $E_{i}$ a uniformly convex Banach
spaces, and $F_{i}:T\to 2^{E_{i}}$ uniformly continuous
set-valued mappings with modulus of continuity $\omega$, i.e.
$$
h((F_{1}(t_{1}),F_{2}(t_{2})),(F_{1}(t_{1}),F_{2}(t_{2})))\le
\omega (\rho (t_{1},t_{2})),\qquad \forall t_{1},t_{2}\in T, \
i=1,2.
$$
Suppose that
images $F_{i}(t)$ are
uniformly convex sets with modulus
$\delta (\ep)$, $i=1,2$ and $\ep\in (0,2r_{0}]$;
$\Delta_{0}=\delta(2r_{0})$. Let $L:E_{1}\oplus E_{2}\to E$ be a
continuous linear surjection and let $\ker L=\L$.

Suppose that
there exists
$C>0$ such that for any $w_{1}=(u_{1},v_{1})\in \L$,
$w_{2}=(u_{2},v_{2})\in \L$ we have $\| u_{1}-u_{2}\|_{E_{1}}\ge
C\| w_{1}-w_{2}\|$ and  $\| v_{1}-v_{2}\|_{E_{2}}\ge C\|
w_{1}-w_{2}\|$.

Let $f(t)\in L(F_{1}(t),F_{2}(t))$ be a uniformly continuous
selection. Then there exists uniformly continuous selections
$f_{i}(t)\in F_{i}(t)$, $i=1,2$, with $f(t)=L(f_{1}(t),f_{2}(t))$.
\end{Th}

\proof The space $E_{1}\oplus E_{2}$ is uniformly convex with the
norm \cite{jap}:
$$\|\cdot \|_{uc}=\sqrt{
\|\cdot\|^{2}_{E_{1}}+\|\cdot\|^{2}_{E_{2}}}.$$ By the
inequalities
$$
\max\{ \| u\|_{E_{1}}, \|v\|_{E_{2}} \}\le\| (u,v) \|_{uc} \le
\sqrt{2}\max\{ \| u\|_{E_{1}}, \|v\|_{E_{2}} \},\qquad\forall
u\in E_{1},\ \forall v\in E_{2}
$$
the norms $\|\cdot \|_{uc}$ and $\max\{ \| u\|_{E_{1}},
\|v\|_{E_{2}} \}$ are equivalent.
 Let $w(t)$ be the metric projection of zero onto
$L^{-1}(f(t))$ in the space $E_{1}\oplus E_{2}$ with the
norm
$\|\cdot \|_{uc}$.

By \cite[Corollary 3.3.6 ]{Aubin} the
set-valued
mapping $t\to L^{-1}(f(t))$ is uniformly continuous
with respect to the Hausdorff distance. By Example 3.2, $w(t)=(u(t),v(t))$ is
uniformly continuous and $L^{-1}(f(t))=w(t)+\L$. Now,
$$
H(t) = w(t)+(F_{1}(t)-u(t),F_{2}(t)-v(t))\cap \L
$$
is uniformly continuous by Lemma 3.\ref{sec}.

We define $(f_{1}(t),f_{2}(t))$ as the metric projection of the
zero onto $H(t)$ in the sense of the norm $\| \cdot \|_{uc}$. This
projection is uniformly continuous by Example 3.2.\qed

\begin{Remark}\label{Best_InRn}
Note that in the case $E=\R^{n}$ we can define
$(f_{1}(t),f_{2}(t))$ as $(f_{1}(t),f_{2}(t))=s(H(t))$, where
$s(H(t))$ is the Steiner point of the set $H(t)$.

Consider set-valued mappings $F_{i}$, $i=1,2$, and the surjection $L$
from Remark 3.\ref{best}, assuming that $\H=\R^{n}$. Let
$f(A)=0\in L(F_{1}(A),F_{2}(A))$. The only solution of this
splitting problem is the point:
$$
f_{1}(A)=f_{2}(A)=p(A)=F_{1}(A)\cap F_{2}(A),
$$
which is the metric projection of zero on the set $A$ in the space
$\R^{n}$. It follows by Remark 3.\ref{best} that in $\R^{n}$ the
order of modulus of continuity for $f_{1}(A)=p(A)$ and
$f_{2}(A)=p(A)$ is exact.

\end{Remark}

\bigskip

\section*{Acknowledgements}
This research was supported by SRA grants P1-0292-0101,
J1-9643-0101, J1-2057-0101, and BI-RU/08-09/001. The first author
was also supported by  RFBR grant 07-01-00156,
ADAP project
"Development of scientific potential of higher school" 2.1.1/500
and  Russian Education Agency project
NK-13P/4. We thank the
referees for comments and suggestions.

\end{document}